\documentclass[conference]{IEEEtran}
\IEEEoverridecommandlockouts
\usepackage{cite}
\usepackage{amsmath,amssymb,amsfonts}
\usepackage{amsthm}
\usepackage{amscd}
\usepackage{accents}
\usepackage{bm}
\usepackage{mathtools}
\usepackage{graphicx}
\usepackage{textcomp}
\usepackage{lipsum}
\usepackage{xcolor}
\usepackage{algorithmicx}
\usepackage{algorithm}
\usepackage{algpseudocode}
\usepackage{pgfplotstable}
\usepackage{setspace}
\usepackage{float}
\usepackage{placeins}
\usepackage{longtable}
\pgfplotsset{compat=1.14}
\definecolor{palette1}{RGB}{228,26,28}
\definecolor{palette2}{RGB}{55,126,184}
\definecolor{palette3}{RGB}{77,175,74}
\definecolor{palette4}{RGB}{152,78,163}
\definecolor{palette5}{RGB}{255,127,0}
\usepackage{hyperref}
\hypersetup{
	colorlinks=true, 
	breaklinks=true,
	urlcolor= blue, 
	linkcolor= blue, 
	citecolor=red, 
	pdftitle={},
	pdfauthor={},
}
\usetikzlibrary{patterns}
\def\BibTeX{{\rm B\kern-.05em{\sc i\kern-.025em b}\kern-.08em
    T\kern-.1667em\lower.7ex\hbox{E}\kern-.125emX}}

\newcommand\ubar[1]{\underaccent{\bar}{#1}}
\newfloat{model}{thp}{lop}
\floatname{model}{Model}

\makeatletter
\let\oldlt\longtable
\let\endoldlt\endlongtable
\def\longtable{\@ifnextchar[\longtable@i \longtable@ii}
\def\longtable@i[#1]{\begin{figure}[t]
\onecolumn
\begin{minipage}{0.5\textwidth}
\oldlt[#1]
}
\def\longtable@ii{\begin{figure}[t]
\onecolumn
\begin{minipage}{0.5\textwidth}
\oldlt
}
\def\endlongtable{\endoldlt
\end{minipage}
\twocolumn
\end{figure}}
\makeatother

\begin{document}

\title{Benchmarking Large-Scale ACOPF Solutions and Optimality Bounds
}

\author{\IEEEauthorblockN{Smitha Gopinath}
\IEEEauthorblockA{\textit{Theoretical Division} \\
\textit{Los Alamos National Laboratory}\\
Los Alamos, USA \\
smitha@lanl.gov}
\and
\IEEEauthorblockN{Hassan L. Hijazi}
\IEEEauthorblockA{\textit{Theoretical Division} \\
\textit{Los Alamos National Laboratory}\\
Los Alamos, USA \\
hlh@lanl.gov}
}

\maketitle

\begin{abstract}
We present the results of a comprehensive benchmarking effort aimed at evaluating and comparing state-of-the-art open-source tools for solving the Alternating-Current Optimal Power Flow (ACOPF) problem.
Our numerical experiments include all instances found in the public library PGLIB with network sizes up to 30,000 nodes. The benchmarked tools span a number of programming languages (Python, Julia, Matlab/Octave, and C$++$), nonlinear optimization solvers (Ipopt, MIPS, and INLP) as well as different mathematical modeling tools (JuMP and Gravity). We also present state-of-the-art optimality bounds obtained using sparsity-exploiting semidefinite programming approaches and corresponding computational times.
\end{abstract}


\section{Introduction}



With the revenues of the energy sector reaching close to $\$400$ billion per year \cite{eia}, power systems optimization plays a central role in maximizing system reliability while minimizing costs.
The Alternating-Current Optimal Power Flow (ACOPF) problem lies at the heart of all power systems optimization. Very recently, the Advanced Research Projects Agency-Energy (ARPA-E) hosted the \href{https://gocompetition.energy.gov}{Grid Optimization Competition} putting the ACOPF at the forefront of impactful research in this field with $\$2.4$ million in prizes at stake.

 
 In this paper, we address the performance of state-of-the-art techniques in providing feasible solutions along with optimality bounds for the ACOPF problem. 
 We compare four open-source implementations capable of producing feasible solutions fast: Gravity (C$++$), PowerModels.jl (Julia), Matpower (Matlab/Octave) and GridOpt (Python). To get lower bounds and compute optimality gaps, we revisit the polynomial semidefinite programming (SDP)- Reformulation-Linearization Technique (RLT) relaxation for ACOPF developed in \cite{Proving}.  While the relaxation has so far been used to prove optimality on all networks up to 300 buses, we present new results on networks up to 30,000 buses. To the best of our knowledge, this is the first published global optimality proof for a number of these instances. We benchmark the performance of this approach against a recent SDP-based moment relaxation introduced by Wang et al. in~\cite{wang2021certifying}.
\section{ACOPF}
\label{sec:ACOPF}
\subsection{Notations} 
\label{sec:notation}
\resizebox{0.5\textwidth}{!}{%
\begin{tabular}{ll}
{\bf Grid Parameters:}\\
$\bm S^g_i = (\bm p^g_i, \bm q^g_i)$  	&Complex, active and reactive power generation at node $i$\\
$\bm S^d_{i} = (\bm p^d_i, \bm q^d_i)$  	&Complex, active and reactive power demand at node $i$ \\
$\bm c_{0i}, \bm c_{1i}, \bm c_{2i}$  			&Generation cost coefficients at node $i$ \\
$\bm t_{ij}$  				&Thermal limit along line $(i,j)$ \\
$\bm Y_{ij}$  				& Complex admittance of line $(i,j)$ $\bm Y_{ij}=\bm g_{ij}+ \bm i \bm b_{ij}$ \\
$\bm Y^c_{ij}$  				&Line charging admittance of line $(i,j)$  \\
$\bm Y^s_{i}$  				&Bus shunt admittance of bus $i$  \\
$\bm T_{ij}$  				&Transformer properties (tap change and phase shift)\\
$\ubar{\bm S}^g_i, \bar{\bm S}^g_i$  	&$\ubar{\bm p}^g_i+\bm i \ubar{\bm q}^g_i, \bar{\bm p}^g_i+\bm i \bar{\bm q}^g_i $ \\
$\ubar{\bm v}_i, \bar{\bm v}_i$  	&Voltage magnitude bounds at node $i$\\
$\ubar{\bm \theta}_{ij} \ge -\pi/2, \bar{\bm \theta}_{ij} \le \pi/2$  &Phase angle difference bounds along line $(i,j)$\vspace{0.2cm}\\
{\bf Grid Variables:}\\
$S_{ij}= ( p_{ij}, q_{ij})$  	&Complex, active and reactive power flow from $i$ to $j$ \\
$V_{i}$  				&Complex voltage at node $i$\\
$I_{ij}$  				&Complex current flowing from $i$ to $j$ at end $i$\\
$l_{ij}$  				&Squared magnitude of current flowing from $i$ to $j$ at end $i$\\
$W_{ij}$    				& Lifted voltage product $V_iV_j^*$ \\
{\bf Other notations:}\\
$\bm {\mathcal G} = (N.E)$ 				&Graph with nodes $N$ and edges $E$\\
$\mathcal{R}(.)$     		& Real component of complex number\\
$\mathcal{I}(.)$     		& Imaginary component of complex number\\
${V}^*$     		& Conjugate of complex number\\
$|V|^2$     		& Square Magnitude of complex number $|V|^2 = \mathcal{R}(V)^2 + \mathcal{I}(V)^2$\\
$\ubar{\bm x}$ , $\bar{\bm x}$       &Lower and upper bounds on $x$\\
\end{tabular}
}
\subsection{Formulations} \label{sec:formulation}

The extended (non-convex) ACOPF problem (with line charging, bus shunts, transformers) is presented in Model \ref{model:ext_ac_opf}.
Note that complex inequalities correspond to component-wise constraints on the respective real and imaginary parts. 

 \begin{model}[!h]
	\caption{Extended ACOPF}
	\label{model:ext_ac_opf}
	{\scriptsize
	\begin{subequations}
		\vspace{-0.2cm}
		\begin{align}
		& \mbox{\bf minimize: } \nonumber \\ 
		& \sum\limits_{i \in N} \bm c_{0i} + \bm c_{1i}\mathcal{R}(S^g_i) + \left(\bm c_{2i}\mathcal{R}(S^g_i)\right)^2 \nonumber\\
		& \mbox{\bf subject to: } \nonumber  \\
		& S_{ij} = (\bm Y^*_{ij}+ \bm Y^{\bm c^*}_{ij})\frac{| V_i |^2}{| \bm T_{ij}|^2}-\bm Y^*_{ij}\frac{V_iV^*_j}{\bm T_{ij}}  ,\quad \forall (i,j) \in E \label{eq:expsij} \\
			& S_{ji} = (\bm Y^*_{ij}+\bm Y^{c^*}_{ij})| V_j |^2-\bm Y^*_{ij}\frac{V_i^*V_j}{\bm T^*_{ij}} ,\quad \forall (i,j) \in E \label{eq:expsji}\\
            & S_i^g - \bm{S}_i^d -\bm{Y^s}_i| V_i |^2= \displaystyle\sum_{(i,j),(j,i) \in E} S_{ij} ,\quad \forall i \in N, \label{eq:exkcls}\\
            &\ubar{\bm v}_i \leq | V_i | \leq \bar{\bm v}_i ,\quad \forall i \in N, \label{eq:exvbounds}\\
		    &\ubar{\bm \theta}_{ij} \mathcal{I}(V_iV^*_j) \leq \mathcal{R}(V_iV^*_j) \leq \bar{\bm \theta}_{ij} \mathcal{I}(V_iV^*_j) ,\quad \forall (i,j) \in E, \label{eq:exanglediff} \\
            &\ubar{\bm S}_i^g \leq S^g_i \leq \bar{\bm S}_i^g ,\quad \forall i \in N, \label{eq:exsbounds} \\
            &| S_{ij} |^2 \leq \bm t_{ij} ,\quad \forall (i,j), (j,i) \in E. \label{eq:exthermal}
		\end{align}
	\end{subequations}
	}
\end{model}

We next present the convex relaxation of Model \ref{model:ext_ac_opf}, developed in \cite{Proving}. In Model \ref{model:sdp_opf}, we introduce the variable $W_{ij}$ to denote the lifted variable $V_iV_j^*$ and $W$ the corresponding matrix. The real component of power flow $S_{ij}$ is denoted by $p_{ij}$, whereas $q_{ij}$ denotes the imaginary component.  Variables $\widehat{p}_{ij}$ and $\widehat{q}_{ij}$ represent squares of real and imaginary power flows on each edge $(i,j)$. Variables $l_{ij}$ and $l_{ji}$ represent the magnitude of current at ends $i$ and $j$, respectively of edge $(i,j)$. Variables $\widehat{lW}_{ij}$, $\widehat{lW}_{ji}$ are introduced to represent products of current magnitude and the squared voltage at each end of an edge.

In \eqref{eq:psd_poly}, $T({\mathcal G})$ is a tree decomposition of a chordal completion of the network graph ${\mathcal G}$ \cite{SDPchordalgraphs}. For each clique (fully-connected sub-graph) $c$ in $T({\mathcal G})$, $W_c$, a submatrix of $W$, contains elements of $W$ that correspond to the nodes in $c$. $M(W_c)$ is the set of all the 2 $\times$ 2 and 3 $\times$ 3 principal minors of $W_c$. Sets ${\mathcal Sec}(y, x)$ and ${\mathcal MC}(x,y,z)$ are defined in \cite{Proving}.

\begin{spacing}{1.15}
Power flow constraints, KCL and voltage-angle constraints are given by \eqref{eq:sdp_pf_from}, \eqref{eq:sdp_pf_to}, \eqref{eq:sdpkcls} and \eqref{eq:sdp_angle}. The thermal limit constraint in lifted variables is given by \eqref{eq:exthermalconvex}. RLT-based valid inequalities on current magnitude are given by \eqref{eq:ext_lij}, \eqref{eq:ext_lji}, \eqref{eq:ext_lWij}, \eqref{eq:ext_lWji}. Polynomial determinant cuts are given in \eqref{eq:psd_poly}. Bound constraints are given in \eqref{eq:exbounds} - \eqref{eq:expbounds} and concave and convex envelopes for  $\widehat{lW}_{ij}$, $\widehat{lW}_{ji}$, $\widehat{p}_{ij}$, $\widehat{q}_{ij}$,  are given in \eqref{eq:ch1} and \eqref{eq:ch2}.
\end{spacing}
\begin{model}[!h]
	\caption{DSDP-RLT}
	\label{model:sdp_opf}
	{\scriptsize
	\begin{subequations}
		\vspace{-0.2cm}
		\begin{align}
		& \mbox{\bf minimize: } \nonumber \\ 
		& \sum\limits_{i \in N} \bm c_{0i} + \bm c_{1i}\mathcal{R}(S^g_i) + \left(\bm c_{2i}\mathcal{R}(S^g_i)\right)^2 \nonumber\\
		& \mbox{\bf subject to: } \nonumber  \\
		& S_{ij} = (\bm Y^*_{ij}+ \bm Y^{\bm c^*}_{ij})\frac{W_{ii}}{| \bm T_{ij}|^2}-\bm Y^*_{ij}\frac{W_{ij}}{\bm T_{ij}}  ,\quad \forall (i,j) \in E \label{eq:sdp_pf_from} \\
		& S_{ji} = (\bm Y^*_{ij}+\bm Y^{c^*}_{ij})W_{jj}-\bm Y^*_{ij}\frac{W_{ij}^*}{\bm T^*_{ij}} ,\quad \forall (i,j) \in E \label{eq:sdp_pf_to}\\
        & S_i^g - \bm{S}_i^d -\bm{Y^s}_i W_{ii} = \displaystyle\sum_{(i,j),  (j,i) \in E} S_{ij} ,\quad \forall i \in N, \label{eq:sdpkcls}\\
		& \mathcal{I}(W_{ij})\tan \ubar{\bm \theta}_{ij} \le \mathcal{R}(W_{ij}) \le \mathcal{I}(W_{ij}) \tan \bar{\bm \theta}_{ij},\quad \forall{(i,j) \in {E}}, \label{eq:sdp_angle}\\
        &\widehat{p_{ij}} + \widehat{q_{ij}} \leq \bm t_{ij}, \quad \forall (i,j), (j,i) \in E. \label{eq:exthermalconvex}\\
 &l_{ij}| \bm T_{ij}| ^2=|(\bm Y_{ij}+\bm{Y^c}_{ij})|^2 W_{ii} -\bm{Y}^*_{ij}(\bm{Y}_{ij}+\bm{Y^c}_{ij})\bm{T}^*_{ij}W_{ij}\nonumber \\
&\;-\bm{Y}_{ij}(\bm{Y}_{ij} +\bm{Y^c}_{ij})^*\bm{T}_{ij}W_{ij}^*+| \bm{Y}_{ij} \bm{T}_{ij}| ^2W_{jj},\quad \forall (i,j) \in E \label{eq:ext_lij}\\
  &\widehat{lW}_{ij} =  | \bm T_{ij}|^2 (\widehat{p}_{ij}+\widehat{q}_{ij}), \quad \forall (i,j) \in E \label{eq:ext_lWij}\\
 & l_{ji}| \bm T_{ij}|^2=|(\bm Y_{ij}+\bm{Y^c}_{ij}) \bm T_{ij}|^2 W_{jj}-\bm Y_{ij}(\bm Y_{ij}+\bm{Y^c}_{ij})^*\bm T_{ij}^*W_{ij} \nonumber \\
&-\bm Y^*_{ij}(\bm Y_{ij}+\bm{Y^c}_{ij})\bm{T}_{ij}W_{ij}^*+| \bm{Y}_{ij}|^2 W_{ii} \forall (i,j) \in E, \label{eq:ext_lji}\\
 & \widehat{lW}_{ji} = \widehat{p}_{ji}+\widehat{q}_{ji} ,\quad \forall (i,j) \in E, \label{eq:ext_lWji}\\
        & \textrm{det}(W_s) \geq 0 \quad \forall W_s \in M(W_c), \forall c \in T({\mathcal G}),  \label{eq:psd_poly}\\
         &\ubar{\bm S}_i^g \leq S^g_i \leq \bar{\bm S}_i^g ,\quad \forall i \in N, \label{eq:exbounds} \\
                 &\ubar{\bm v}^2_i \leq |W_{ii}| \leq \bar{\bm v}^2_i ,\quad \forall i \in N, \label{eq:sdpvbounds}\\
        &\bar{\bm v}_i\bar{\bm v}_j \cos \ubar{\bm \theta}_{ij} \leq \mathcal{R}(W_{ij}) \leq \bar{\bm v}_i\bar{\bm v}_j ,\quad \forall (i,j) \in E, \label{eq:sdpRwijbounds}\\
        &-\bar{\bm v}_i\bar{\bm v}_j \leq \mathcal{I}(W_{ij}) \leq \bar{\bm v}_i\bar{\bm v}_j ,\quad \forall (i,j) \in E, \label{eq:sdpIwijbounds}\\
        & 0 \le l_{ij} \le |\bm T_{ij}|^2\bm{t}_{ij}/\ubar{\bm v}^2_i \quad \forall (i,j) \in E\\
         & 0 \le l_{ji} \le \bm{t}_{ij}/\ubar{\bm v}^2_j \quad \forall (i,j) \in E\\
          & -\sqrt{\bm{t}_{ij}} \le p_{ij}, q_{ij} \le \sqrt{\bm{t}_{ij}} \quad \forall (i,j), (j,i) \in E \label{eq:expbounds}\\
        & {\mathcal MC}(l_{ij},W_{ii},\widehat{lW}_{ij}); {\mathcal MC}(l_{ji},W_{jj},\widehat{lW}_{ji}) ,\quad \forall (i,j) \in E \label{eq:ch1},\\
        & {\mathcal Sec}(\widehat{p}_{ij}, p_{ij}),{\mathcal Sec}(\widehat{q}_{ij}, q_{ij}), \quad \forall (i,j), (j,i) \in E \label{eq:ch2}
		\end{align}
	\end{subequations}
	}
\end{model}

\subsection{Scaling Up the Determinant-Based Relaxations}\label{sec:scaling}
Below are the steps we took to make sure Model \eqref{model:sdp_opf} can still converge on large-scale instances:
\begin{enumerate}
\item The total number of determinant cuts \eqref{eq:psd_poly} generated from 3 $\times$ 3 determinants are capped at 30,000.
\item The RLT constraints on current-magnitudes \eqref{eq:ext_lij},\eqref{eq:ext_lWij}, \eqref{eq:ext_lji}, \eqref{eq:ext_lWji}  are applied on the edge set $E$' where:\\
$E'=\{(i,j) \in E: |\bm{g}_{ij}| \le 10^{3}, |\bm{b}_{ij}| \le 10^{3}\}$
\item Model \eqref{model:sdp_opf} is initialized using the solution of Model \eqref{model:ext_ac_opf} and the non-convex relationships between the variables.
\item Concave envelopes for $\widehat{p}_{ij}, \, \widehat{q}_{ij} \,\, \forall (i,j), (j,i) \in E$ were not added to the formulation (secant constraints in \cite{Proving}).
\item Thermal limit constraints \eqref{eq:exthermal}, \eqref{eq:exthermalconvex} and the objective function are scaled by $10^{-3}$.
\end{enumerate}

\section{Computational Experiments}
Upper bounds were computed for all instances in PGLiB v21.07 \cite{pglib}. We tested Gravity \cite{Gravity} and PowerModels v0.18.3 \cite{PowerModels} with Ipopt 3.12~\cite{Ipopt} and linear solver MA27 \cite{HSL}, GridOpt v1.3.7 \cite{GridOpt} with the default interior-point solver INLP from the python-package OPTALG \cite{INLP}, and MATPOWER v7.1 \cite{Matpower} which uses the default interior-point solver MIPS. 
We used a Mac 3.5 GHz 6-Core Intel Xeon E5 with 32 GB of RAM for computing all ACOPF solutions. 

Lower bounds for all networks in PGLiB v21.07, with a second order cone relaxation gap greater than 1$\%$, were tested. Lower bounds from Model \eqref{model:sdp_opf} were computed in Gravity~\cite{Gravity}, using Ipopt 3.12 \cite{Ipopt} with linear solver MA57 \cite{HSL}.  The optimality and constraint satisfaction tolerances are set to $10^{-6}$ and bound relax factor to $10^{-8}$. Times reported include lower bound model build time, initialization time and solution time. Lower bounds corresponding to the moment-hierarchy TSSOS \cite{TSSOS} were tested in Julia 1.6 using Mosek 9.3 for solving the underlying SDPs. Times include model build time, prepossessing and solution time. Lower bounds are run on a 3.5 GHz Dual-Core Intel Core i7 with 16 GB RAM. In all experiments, we exclude Julia's just-in-time compile time.
All experiments presented here can be reproduced by cloning the following repositories: \href{https://github.com/coin-or/Gravity/tree/PES22}{PES22} and \href{https://github.com/coin-or/Gravity/tree/PES22_Largescale}{{PES22\_Largescale}}.

\section{Numerical Results}
\label{Results}
\subsection{Upper Bounds}
\begin{figure*}[!htbp]
\centering
\resizebox{0.6\textwidth}{!}{
\begin{tikzpicture}[baseline]
    \begin{axis}[ybar = .05cm,
    bar width = 1.5pt,
    enlarge x limits = {value = .25, upper},
    xmin = 0,
    xmax = 60,
    ymin = 0,
    ymax = 15,enlargelimits=false,
	                     xtick pos=left,ytick pos=left,height=10.0cm,
	                            ylabel=$\times$ times slower than Gravity,
								xlabel={Instance number},
								legend image post style={scale=0.5},
								legend style={font=\small,at={(0.3,0.95)},anchor=north},
								legend cell align={left},
								legend columns=3
    ]
        \pgfplotstableread[col sep = space]{CSVs/Matpower_Gravity.csv}{\MGravity};
        \addplot[palette1,fill=palette1] table [x = Num, y = Factor]{\MGravity};
        \addlegendentry{Matpower}
        \pgfplotstableread[col sep = space]{CSVs/PowerModels_Gravity.csv}{\PMGravity};
        \addplot[palette2, semithick] table [x = Num, y = Factor]{\PMGravity};
        \addlegendentry{PowerModels}
        \pgfplotstableread[col sep = space]{CSVs/GridOpt_Gravity.csv}{\GOGravity};
        \addplot[palette3, pattern=horizontal lines] table [x = Num, y = Factor]{\GOGravity};
        \addlegendentry{GridOpt}
    \end{axis}
\end{tikzpicture}\hfill
}
	 \caption{Relative Slowdown Factor Compared to Gravity on ACOPF Typical Operating Conditions (TYP) PGLIB instances.}
    \label{fig:Relative}
\end{figure*}
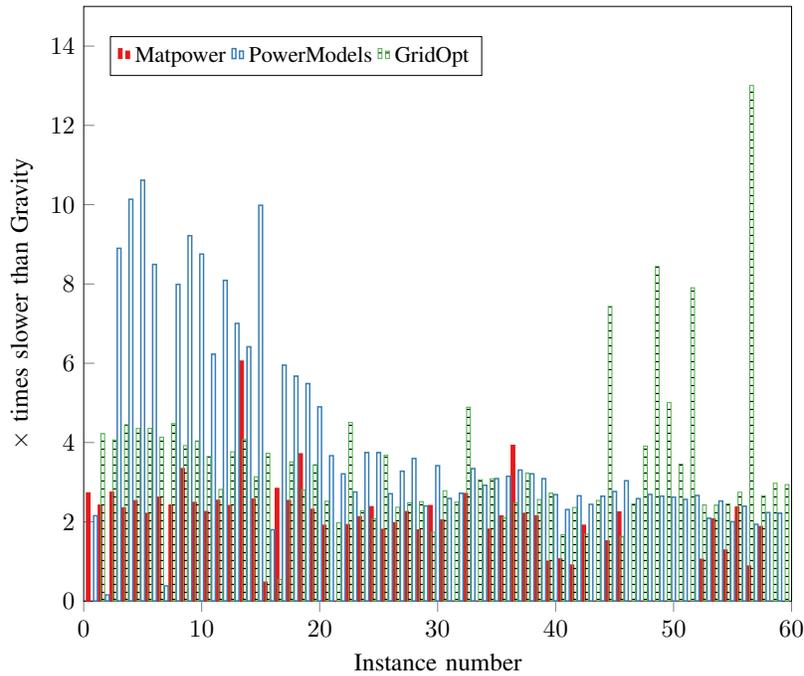

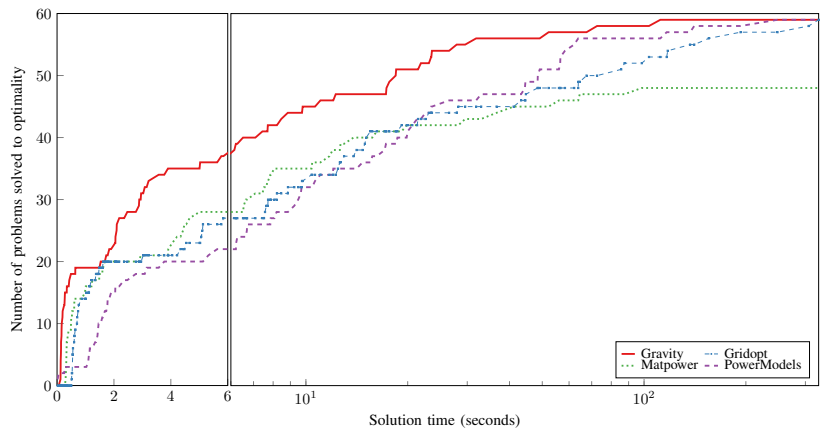
\begin{figure*}[!htbp]
\centering
\resizebox{0.6\textwidth}{!}{
\begin{tikzpicture}[baseline]
    \begin{axis}[legend cell align=left,enlargelimits=false,xtick pos=left,
	              ytick pos=left,height=10.0cm,ylabel=Number of problems solved to optimality,
                 width=0.3\textwidth,,ymin=0,ymax=60,restrict x to domain=0:350,xmin=0.0,xmax=6]
        \pgfplotstableread[col sep = space]{CSVs/Gravity_PP_TYP.csv}{\Gravity};
        \addplot[very thick, palette1] table [x =solve_time, y =num_problems_solved]{\Gravity};
        \pgfplotstableread[col sep = space]{CSVs/Gridopt_PP_TYP.csv}{\Gridopt};
        \addplot[dashed,mark=square*,mark options={scale=0.2,solid}, palette2] table [x = solve_time, y = num_problems_solved]{\Gridopt};
        \pgfplotstableread[col sep = space]{CSVs/Matpower_PP_TYP.csv}{\Matpower};
        \addplot[dotted, very thick, palette3] table [x = solve_time, y = num_problems_solved]{\Matpower};
        \pgfplotstableread[col sep = space]{CSVs/PowerModels_PP_TYP.csv}{\PowerModels};
        \addplot[dashed, very thick, palette4] table [x = solve_time, y = num_problems_solved]{\PowerModels};
    \end{axis}
\end{tikzpicture}\hspace{-0.75em}
\begin{tikzpicture}[baseline]
	\begin{semilogxaxis}[legend cell align=left,enlargelimits=false,
	                     xtick pos=left,ytick pos=left,height=10.0cm,
								xlabel={\hspace{-0.2\textwidth}Solution time (seconds)},
								ytick=\empty,width=0.82\textwidth,
								xlabel absolute,legend style={at={(0.985, 0.02)},anchor=south east},
								legend image post style={scale=0.5},
								legend style={font=\small,column sep=1.5pt, row sep=-4.0pt},
								legend columns=2,ymin=0,ymax=60,xmin=6.0]
        \pgfplotstableread[col sep = space]{CSVs/Gravity_PP_TYP.csv}{\Gravity};
        \addplot[very thick, palette1] table [x = solve_time, y = num_problems_solved]{\Gravity};
        \addlegendentry{Gravity}
        \pgfplotstableread[col sep = space]{CSVs/Gridopt_PP_TYP.csv}{\Gridopt};
        \addplot[dashed,mark=square*,mark options={scale=0.2,solid},palette2] table [x = solve_time, y = num_problems_solved]{\Gridopt};
        \addlegendentry{Gridopt}
        \pgfplotstableread[col sep = space]{CSVs/Matpower_PP_TYP.csv}{\Matpower};
        \addplot[dotted,very thick, palette3] table [x = solve_time, y = num_problems_solved]{\Matpower};
        \addlegendentry{Matpower}
        \pgfplotstableread[col sep = space]{CSVs/PowerModels_PP_TYP.csv}{\PowerModels};
        \addplot[dashed, very thick, palette4] table [x = solve_time, y = num_problems_solved]{\PowerModels};
        \addlegendentry{PowerModels}
    \end{semilogxaxis}
\end{tikzpicture}\hfill
}
	 \caption{Performance Profile on Typical Operating Conditions (TYP) PGLIB Instances.}
    \label{fig:acopf_typ}
\end{figure*}

\begin{figure}[!htbp]
\centering
\resizebox{0.5\textwidth}{!}{
\begin{tikzpicture}[baseline]
    \begin{axis}[ybar = .05cm,
    bar width = 1.5pt,
    enlarge x limits = {value = .25, upper},
    xmin = 0,
    xmax = 60,
    ymin = -100,
    ymax = 600,enlargelimits=false,
	                     xtick pos=left,ytick pos=left,height=10.0cm,
	                            ylabel=$\%$ slower than scaled version,
								xlabel={Instance number},
								legend image post style={scale=0.5},
								legend style={font=\small,at={(0.8,0.95)},anchor=north},
								legend cell align={left},
								legend columns=1
    ]
        \pgfplotstableread[col sep = space]{CSVs/Scaling.csv}{\MGravity};
        \addplot[palette1, semithick] table [x = Num, y = Factor]{\MGravity};
        \addlegendentry{Scaling Off}
    \end{axis}
\end{tikzpicture}\hfill
}
	 \caption{Scaling the Objective and Thermal Limits Constraints by $10^{-3}$.}
    \label{fig:Scaling}
\end{figure}
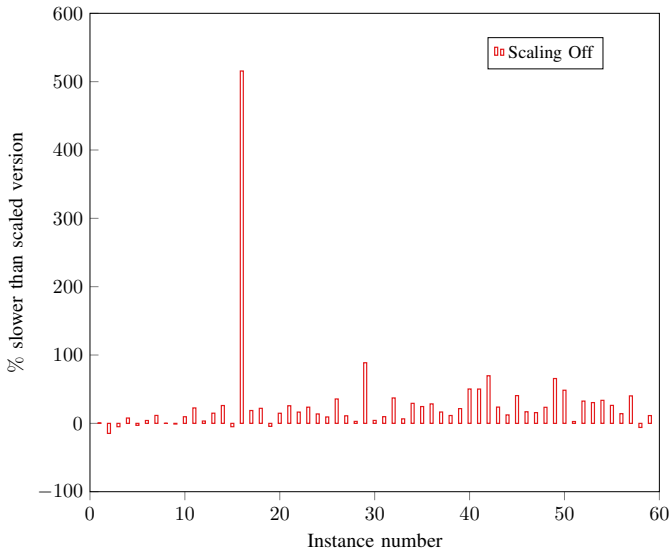

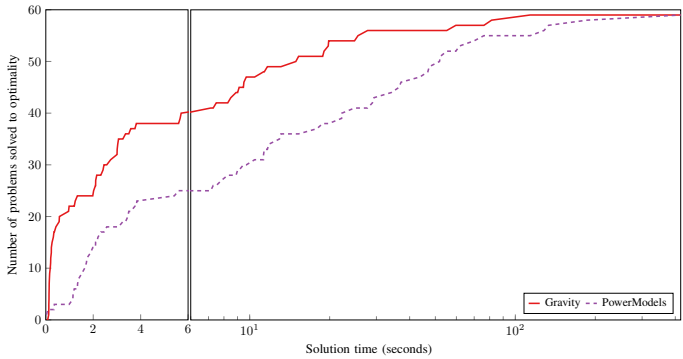
\begin{figure}[!htbp]
\centering
\resizebox{0.5\textwidth}{!}{
\begin{tikzpicture}[baseline]
    \begin{axis}[legend cell align=left,enlargelimits=false,xtick pos=left,
	              ytick pos=left,height=10.0cm,ylabel=Number of problems solved to optimality,
                 width=0.3\textwidth,,ymin=0,ymax=60,restrict x to domain=0:420,xmin=0.0,xmax=6]
        \pgfplotstableread[col sep = space]{CSVs/Gravity_PP_API.csv}{\Gravity};
        \addplot[very thick, palette1] table [x = solve_time, y = num_problems_solved]{\Gravity};
        \pgfplotstableread[col sep = space]{CSVs/PowerModels_PP_API.csv}{\PowerModels};
        \addplot[dashed, very thick, palette4] table [x = solve_time, y = num_problems_solved]{\PowerModels};
    \end{axis}
\end{tikzpicture}\hspace{-0.75em}
\begin{tikzpicture}[baseline]
	\begin{semilogxaxis}[legend cell align=left,enlargelimits=false,
	                     xtick pos=left,ytick pos=left,height=10.0cm,
								xlabel={\hspace{-0.2\textwidth}Solution time (seconds)},
								ytick=\empty,width=0.82\textwidth,
								xlabel absolute,legend style={at={(0.985, 0.02)},anchor=south east},
								legend image post style={scale=0.5},
								legend style={font=\small,column sep=1.5pt, row sep=-4.0pt},
								legend columns=2,ymin=0,ymax=60,xmin=6.0]
        \pgfplotstableread[col sep = space]{CSVs/Gravity_PP_API.csv}{\Gravity};
        \addplot[very thick, palette1] table [x = solve_time, y = num_problems_solved]{\Gravity};
        \addlegendentry{Gravity}
        \pgfplotstableread[col sep = space]{CSVs/PowerModels_PP_API.csv}{\PowerModels};
        \addplot[dashed, very thick, palette4] table [x = solve_time, y = num_problems_solved]{\PowerModels};
        \addlegendentry{PowerModels}
    \end{semilogxaxis}
\end{tikzpicture}\hfill
}
	 \caption{Performance Profile on Congested (API) PGLIB Instances.}
    \label{fig:acopf_api}
\end{figure}

\begin{figure}[!htbp]
\centering
\resizebox{0.5\textwidth}{!}{
\begin{tikzpicture}[baseline]
    \begin{axis}[legend cell align=left,enlargelimits=false,xtick pos=left,
	              ytick pos=left,height=10.0cm,ylabel=Number of problems solved to optimality,
                 width=0.3\textwidth,,ymin=0,ymax=60,restrict x to domain=0:420,xmin=0.0,xmax=6]
        \pgfplotstableread[col sep = space]{CSVs/Gravity_PP_sad.csv}{\Gravity};
        \addplot[very thick, palette1] table [x = solve_time, y = num_problems_solved]{\Gravity};
        \pgfplotstableread[col sep = space]{CSVs/PowerModels_PP_sad.csv}{\PowerModels};
        \addplot[dashed, very thick, palette4] table [x = solve_time, y = num_problems_solved]{\PowerModels};
    \end{axis}
\end{tikzpicture}\hspace{-0.75em}
\begin{tikzpicture}[baseline]
	\begin{semilogxaxis}[legend cell align=left,enlargelimits=false,
	                     xtick pos=left,ytick pos=left,height=10.0cm,
								xlabel={\hspace{-0.2\textwidth}Solution time (seconds)},
								ytick=\empty,width=0.82\textwidth,
								xlabel absolute,legend style={at={(0.985, 0.02)},anchor=south east},
								legend image post style={scale=0.5},
								legend style={font=\small,column sep=1.5pt, row sep=-4.0pt},
								legend columns=2,ymin=0,ymax=60,xmin=6.0]
        \pgfplotstableread[col sep = space]{CSVs/Gravity_PP_sad.csv}{\Gravity};
        \addplot[very thick, palette1] table [x = solve_time, y = num_problems_solved]{\Gravity};
        \addlegendentry{Gravity}
        \pgfplotstableread[col sep = space]{CSVs/PowerModels_PP_sad.csv}{\PowerModels};
        \addplot[dashed, very thick, palette4] table [x = solve_time, y = num_problems_solved]{\PowerModels};
        \addlegendentry{PowerModels}
    \end{semilogxaxis}
\end{tikzpicture}\hfill
}
	 \caption{Performance Profile on Small Angle Difference (SAD) PGLIB Instances.}
    \label{fig:acopf_sad}
\end{figure}
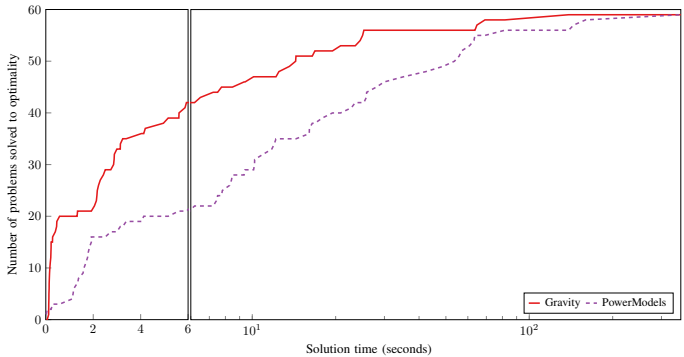

Table \ref{Table:ub} presents the wall-clock computational time for solving the ACOPF problem on the original (typical operating conditions) PGLIB v21.07 library \cite{pglib}.
Figure \ref{fig:Relative} shows the slow-down of all solvers compared to the Gravity-Ipopt approach. The figure shows that the Python-based tool GridOpt can be up to 13$\times$ slower with an average 3.5$\times$ slow-down. Powermodels, using the same Ipopt solver can be up to 10$\times$ slower with an average 4$\times$ slow-down. Finally, note that Matpower fails on 11 out of 59 instances. When it does converge, Matpower is on average 2.2$\times$ slower with a worst-case slow-down of 6$\times$. Performance profiles depicting the number of instances solved as a function of time is presented on all PGLIB instances in Figures \ref{fig:acopf_typ}, \ref{fig:acopf_api} and \ref{fig:acopf_sad}. Note that scaling the objective function and the thermal limit constraints by $10^{-3}$ helps speed-up Ipopt's convergence by $28\%$ on average as illustrated in Figure \ref{fig:Scaling}. All the upper-bound results presented here (except for Figure~\ref{fig:Scaling}) are without scaling. All the semidefinite programming relaxations used to compute lower bounds use this scaling as discussed in Section~\ref{sec:scaling}.

\subsection{Lower Bounds}
Results of the lower bound runs are given in Table \ref{Table:lb}. We compute the optimality gap using the optimal solution of Model \eqref{model:sdp_opf} with and without the RLT cuts. For networks with an optimality gap greater than 1\% and less than a 1000 buses, we also compute the gap from the optimization-based bound tightening (OBBT) algorithm (SDP-BT) using a small cluster of 8 machines, with a time limit of 3h, as described in \cite{Proving}. 

Columns 1 and 2 give the root node gap and computation time in seconds, respectively for the Determinant-SDP (DSDP) relaxation. Columns 3 and 4 give the root node gap and computation time in seconds, respectively for the Determinant-SDP relaxation with the addition of the RLT constraints. If OBBT is called, Columns 5 and 6 give the gap and time after bound-tightening.

F denotes that Ipopt failed to converge.
T.L. denotes that Ipopt failed to obtain a solution within the maximum time of 1.5 h for the relaxation.
$*$ denotes that Mosek reported problem status as feasible (not optimal).
$\#$ denotes that Mosek reported problem status as unknown. MEM denotes insufficient memory.
Gaps in bold indicate instances in which the DSDP approaches close the gap, but the moment hierarchy does not.

For small and medium-scale instances ($\le$ 3000 buses), the DSDP-RLT relaxation closes the gap on 61\% of all tested instances. We ran OBBT on 17 out of the remaining 24 open cases, which closed the gap on all instances except two. Thus, the DSDP-RLT, in conjunction with OBBT, can successfully close the gap on 85\% on all instances up to 3000 buses in PGLiB.

The DSDP relaxation, without the RLT cuts, is more robust in solving large-scale instances ($>$ 3000 buses). Together, the two DSDP relaxations converge to a lower bound for all instances, except one (case6515\_rte\_api). In total, the DSDP and DSDP-RLT relaxations, close the gap for 53 out of the 118 open cases in PGLiB.

We were able to test 57 small to medium-scale PGLiB cases with the moment-hierarchy approach, without running into memory issues on our hardware. For these 57 instances, the DSDP-RLT relaxation and the SDP-BT close the gap for 32 and 49 cases, respectively. In contrast, the level 1 and the two levels combined close the gap for 11 and 33 cases, respectively. However, the solution status for several of these problems is ``feasible", and not certified to be optimal.

{
\begin{table}		
\resizebox{0.5\textwidth}{!}{
\begin{tabular}{|l|r|rrrr|}
\hline
Instance & Obj & Matpower & PowerModels & Gravity &  Gridopt\\ 
\hline
\hline
 case3\_lmbd    &  5812.64 &  0.3 &  0.3 &  0.1 &  0.5\\
 case5\_pjm    &  17551.89 &  0.3 &  0.1 &  0.1 &  0.5\\
 case14\_ieee    &  2178.08 &  0.3 &  1.1 &  0.1 &  0.6\\
 case24\_ieee\_rts    &  63352.2 &  0.3 &  1.4 &  0.1 &  0.6\\
 case30\_as    &  803.13 &  0.3 &  1.3 &  0.1 &  0.6\\
 case30\_ieee    &  8208.52 &  0.3 &  1.2 &  0.1 &  0.6\\
 case39\_epri    &  138415.56 &  0.4 &  0.1 &  0.2 &  0.7\\
 case57\_ieee    &  37589.34 &  0.3 &  1.1 &  0.1 &  0.6\\
 case60\_c    &  92693.67 &  0.5 &  1.5 &  0.2 &  0.6\\
 case73\_ieee\_rts    &  189764.08 &  0.4 &  1.4 &  0.2 &  0.6\\
 case89\_pegase    &  107285.67 &  0.6 &  1.7 &  0.3 &  0.8\\
 case118\_ieee    &  97213.61 &  0.5 &  1.6 &  0.2 &  0.7\\
 case162\_ieee\_dtc    &  108075.64 &  0.6 &  1.8 &  0.3 &  1\\
 case179\_goc    &  754266.41 &  1.7 &  1.8 &  0.3 &  0.9\\
 case200\_activ    &  27557.57 &  0.5 &  1.9 &  0.2 &  0.7\\
 case240\_pserc    &  3329670.06 &  1 &  3.8 &  2.1 &  1.1\\
 case300\_ieee    &  565219.97 &  1 &  2.1 &  0.3 &  1.2\\
 case500\_goc    &  454945.98 &  1.2 &  2.8 &  0.5 &  1.4\\
 case588\_sdet    &  313139.78 &  1.6 &  2.4 &  0.4 &  1.5\\
 case793\_goc    &  260197.85 &  1.5 &  3.2 &  0.6 &  1.6\\
 case1354\_pegase    &  1258843.99 &  2.9 &  5.6 &  1.5 &  3\\
 case1888\_rte    &  1402530.83 &     $\#$   &  9.2 &  2.9 &  13\\
 case1951\_rte    &  2085581.84 &  12 &  17.1 &  6.2 &  14.1\\
 case2000\_goc    &  973432.47 &  4.7 &  8.2 &  2.2 &  4.5\\
 case2312\_goc    &  441330.33 &  5 &  7.9 &  2.1 &  7.7\\
 case2383wp\_k    &  1868191.59 &  4.5 &  6.7 &  2.5 &  5.9\\
 case2736sp\_k    &  1308015 &  4.1 &  6.7 &  2 &  5.1\\
 case2737sop\_k    &  777727.69 &  3.9 &  6.2 &  1.7 &  4.3\\
 case2742\_goc    &  275705.45 &  19.9 &  26.5 &  11 &  19.2\\
 case2746wop\_k    &  1208258.5 &  4.4 &  6.3 &  1.8 &  5.1\\
 case2746wp\_k    &  1631707.93 &  4.2 &  5.3 &  2.1 &  5.1\\
 case2848\_rte    &  1286608.2 &     $\#$   &  15.7 &  5.8 &  28.3\\
 case2853\_sdet    &  2052386.73 &  7.9 &  9.7 &  2.9 &  8.9\\
 case2868\_rte    &  2009605.31 &     $\#$   &  14.7 &  5 &  15.5\\
 case2869\_pegase    &  2462790.44 &  7.1 &  12.1 &  3.9 &  8.2\\
 case3012wp\_k    &  2600842.73 &  6.6 &  9.7 &  3.1 &  7.6\\
 case3022\_goc    &  601383.85 &  12.7 &  10.6 &  3.2 &  10.4\\
 case3120sp\_k    &  2147969.08 &  6.6 &  9.5 &  3 &  7.6\\
 case3375wp\_k    &  7438169.42 &  7.7 &  11 &  3.6 &  9.7\\
 case3970\_goc    &  960985.26 &  7.5 &  20 &  7.4 &  12.4\\
 case4020\_goc    &  822247.29 &  10.4 &  22.6 &  9.8 &  23.1\\
 case4601\_goc    &  826241.53 &  8.1 &  23.5 &  8.8 &  14.9\\
 case4619\_goc    &  476703.73 &  16.2 &  20.6 &  8.5 &  21.5\\
 case4661\_sdet    &  2251344.07 &     $\#$   &  17.3 &  6.5 &  48.4\\
 case4837\_goc    &  872255.32 &  11.7 &  21.3 &  7.7 &  12.6\\
 case4917\_goc    &  1387791.02 &  13.8 &  18.7 &  6.1 &  15.1\\
 case6468\_rte    &  2069730.15 &     $\#$   &  56.6 &  21.9 &  85.6\\
 case6470\_rte    &  2237570.59 &     $\#$   &  33 &  12.3 &  103.4\\
 case6495\_rte    &  3067825.31 &     $\#$   &  62.1 &  23.5 &  117.5\\
 case6515\_rte    &  2825499.65 &     $\#$   &  48.5 &  18.5 &  63.9\\
 case8387\_pegase   &  2771392.29 &     $\#$   &  44.4 &  17.3 &  136.8\\
 case9241\_pegase    &  6243090.4 &     $\#$   &  49.1 &  18.4 &  44.7\\
 case9591\_goc    &  1061683.56 &  29.6 &  58.4 &  27.9 &  67.6\\
 case10000\_goc    &  1354031.34 &  36.7 &  44.6 &  17.7 &  43.3\\
 case10480\_goc    &  2314648.01 &  41.2 &  63.9 &  31.9 &  87.5\\
 case13659\_pegase    &  8948048.97 &  56 &  56.5 &  23.6 &  306.8\\
 case19402\_goc    &  1977815.4 &  64.5 &  140.7 &  72.6 &  193\\
 case24464\_goc    &  2629531.3 &  98.4 &  116.6 &  52.2 &  155.3\\
 case30000\_goc    &  1142331.6 &     $\#$   &  248.1 &  111.7 &  328.3\\

\hline
\end{tabular}
}
\vspace{0.05cm}
\caption{ACOPF PGLIB v21.07 instances. Running time in seconds.\label{Table:ub}}
\vspace{-0.5cm}
\end{table}	
}
{
\begin{table}		
\resizebox{0.49\textwidth}{!}{
\begin{tabular}{|l|rrrrrr|rrrr|}																									
\hline																
	&	\multicolumn{6}{c|}{DSDP-RLT-OBBT}							&	\multicolumn{4}{c|}{Moment Hierarchy}\\

Instance & DSDP	&	time1	&	DSDP-RLT	&	time2	& SDP-BT & time & Gap1	&	time1	&	Gap2	&	time2	 \\
\hline
nesta\_case9\_bgm\_nco & 	10.8 &	0.1 & 	10.8 &	0.1 & 	\textbf{0.1} &	1.4 &	10.8 &	21.9 &	10.8 &	0.8 \\
\hline
case3\_lmbd  & 	0.4 &	0 &	0 &	0.1 &	 &	  & 	1.2 & 	22.5 & 	0.0 & 	0.8 \\ 
case5\_pjm & 	5.2 &	0 &	0.1 &	0.1 &	  & 	  & 	14.6 & 	19.5 & 	0.1 & 	0.8 \\ 
case30\_ieee & 	0.0 &	0.1 &	0 &	0.2 &	   & 	 & 	8.1 & 	23.7 & 	0.0 & 	1.6 \\ 
case162\_ieee\_dtc & 	1.8 &	2.3 &	1.6 &	4 &	\textbf{0.1} & 	2760.2 & 	6.0 & 	21.0 & 	1.5 & 	103.6 \\ 
case240\_pserc & 	1.5 &	1.5 &	1.2 &	8.1 &	0.6 & 	9676.0 & 	2.4$^*$ & 	20.6 & 	0.6$^*$ & 	31 \\ 
case300\_ieee & 	0.2 &	2.4 &	0.1 &	2.3 &	   & 	   & 	1.9$^*$ & 	23.3 & 	0$^*$ & 	29.5 \\ 
case588\_sdet & 	1.0 &	1.9 &	0.4 &	9.9 &	 & 	  & 	1.4$^*$ & 	25.8 & 	0.4$^*$ & 	63 \\ 
case793\_goc & 	0.8 &	20.1 &	0.6 &	5.0 &	 & 	 & 	1.5 & 	25.4 & 	0.3$^*$ & 	85.8 \\ 
case1888\_rte & 	2.0 &	12.4 &	2.0 &	139.8 &	 & 	 & 	2.6$^*$ & 	50.0 & 	1.9$^*$ & 	565.3 \\ 
case2312\_goc & 	1.2 &	77.7 &	1.0 &	128.8 &	 & 	 & 	1.6$^*$ & 	149.3 & 	0.6$^*$ & 	988 \\ 
case2383wp\_k & 	0.3 &	42.6 &	0.3 &	80.8 &	 & 	  & 	0.5$^*$ & 	84.6 & 	 & 	\\ 
case2742\_goc & 	0.2 &	74.4 &	0.2 &	246.3 &	 & 	 & 	0.0$^*$ & 	837.2 & 	 & 	 \\ 
case2869\_pegase & 	\textbf{0.1} &	28.5 &	F &	 &	 & 	 & 	MEM & 	 & 	 & 	\\ 
\hline
case3012wp\_k & 	\textbf{0.3} &	2044.6 &	\textbf{0.2} &	635.5 &	 & 	 & 	 	MEM & 	 & 	 & 	\\
case3022\_goc & 	1.9 &	39.5 &	F &	146.2 &	 & 	 & 		MEM & 	 & 	 & 	\\
case4020\_goc & 	\textbf{0.8} &	119.5 &	\textbf{0.8} &	334 &	 & 	 & 		MEM & 	 & 	 & 	\\
case4661\_sdet & 	0.9 & 914.6	 &	\textbf{0.9} &	822.1 &	 & 	 & 	MEM	 & 	 & 	 & 	\\
case4917\_goc & 	1.5 &	208.6 &	\textbf{0.9} &	175.8 &	 & 	 & 		MEM & 	 & 	 & 	\\
case6468\_rte & 	\textbf{0.6} &	208.5 &	F &	 &	 & 	 & 	MEM	 & 	 & 	 & 	\\
case6470\_rte & 	\textbf{0.9} &	329.9 &	 & T.L.	 &	 & 	 & 	MEM	 & 	 & 	 & 	\\
case6495\_rte & 	13.3 &	192.4 &	 &	T.L. &	 & 	 & 		MEM & 	 & 	 & 	\\
case6515\_rte & 	5.6 &	230.3 &	 &	T.L. &	 & 	 & 		MEM & 	 & 	 & 	\\
case9241\_pegase & 	2.2 &	251.9 &	F &	 &	 & 	 & 	MEM	 & 	 & 	 & 	\\
case10000\_goc & 	1.5 &	667.9 &	1.5 &	272.9 &	 & 	 & 	MEM	 & 	 & 	 & 	\\
case10480\_goc & 	1.1 &	139 &	1.1 &	554.4 &	 & 	 & 	MEM	 & 	 & 	 & 	\\
case13659\_pegase & 	1.2 &	247.2 &	F &	773.7 &	 & 	 & 	MEM	 & 	 & 	 & 	\\
case19402\_goc & 	1.2 &	344.3 &	1.2 &	853 &	 & 	 & 	MEM	 & 	 & 	 & 	\\
case30000\_goc & 	2.9 &	1117.8 &	2.9 &	792.1 &	 & 	 & 	MEM	 & 	 & 	 & 	\\
\hline
\hline
case3\_lmbd\_api & 	7.1 &	0.0 &	0.9 &	0.1 &	  &	 &	4.9 & 	24.6 & 	0.0 & 	0.7 \\ 
case5\_pjm\_api & 	0.3 &	0.0 &	0.1 &	0.1 &	  &	  &	4.1 & 	19.3 & 	0.0 & 	0.8 \\ 
case14\_ieee\_api & 	0.1 &	0.1 &	0.1 &	0.1 &	  &	  &	5.2 & 	18.8 & 	0.0 & 	0.9 \\ 
case24\_ieee\_rts\_api & 	2.2 &	0.0 &	1.1 &	0.1 &	0.0 & 	2.1 &	6.4 & 	19.0 & 	0.8 & 	1.9 \\ 
case30\_as\_api & 	6.9 &	0.1 &	5.8 &	0.2 &	0.7 & 	0.3 & 	42.9 & 	21.6 & 	-0.2$^*$ & 	3.4 \\ 
case30\_ieee\_api & 	0.2 &	0.1 &	0.2 &	0.2 &	  &	  &	4.4 & 	20 & 	0.0 & 	1.8 \\ 
case39\_epri\_api & 	0.2 &	0.1 &	0.2 &	0.2 &	  &	  &	1.8 & 	22.7 & 	0.0 & 	4.8 \\ 
case73\_ieee\_rts\_api & 	3.0 &	0.3 &	2.2 &	0.5 &	\textbf{0.4} & 	25.0 & 	5.6$^*$ & 	20.2 & 	1.8$^*$ & 	7.0 \\ 
case89\_pegase\_api & 	21.9 &	0.6 &	21.7 &	5.2 &	\textbf{0.5} & 	545.2 & 	22.8 & 	20.4 & 	21.7$^*$ & 	1879.3 \\ 
case118\_ieee\_api & 	12.1 &	0.3 &	9.2 &	1.0 &\textbf{	0.9} & 	116.2 & 	20 & 	25.2 & 	8.6 & 	10.9 \\ 
case162\_ieee\_dtc\_api & 	1.5 &	2.8 &	1.3 &	3.8 &	\textbf{0.0} & 	1282.7 & 	7.4$^*$ & 	20.5 & 	1.2$^*$ & 	109.5 \\ 
case179\_goc\_api & 	0.6 &	0.9 &	0.6 &	2.1 &	  & 	  & 	10.8$^*$ & 	19.7 & 	0.5$^*$ & 	10.1 \\ 
case500\_goc\_api & 	2.4 &	3.7 &	2.2 &	6.8 &	\textbf{0.6} & 	1953.8 & 	4.7 & 	27.6 & 	2$^*$ & 	89.7 \\ 
case588\_sdet\_api & 	0.5 &	3.9 &	0.2 &	9.0 &	  & 	  & 	1.1$^*$ & 	24.3 & 	0.2$^*$ & 	62.1 \\ 
case793\_goc\_api & 	6.4 &	9.3 &	2.9 &	3.5 &	\textbf{0.9} & 	3269.4 & 	6.6 & 	27.3 & 	1.6 & 	90.4 \\ 
case2000\_goc\_api & 	0.7 &	199.8 &	0.7 &	129 &	 & 	 & 	3.7$^*$ & 	74.8 & 	0.5$^*$ & 	1204.5 \\ 
case2312\_goc\_api & 	9.9 &	123.3 &	8.6 &	143.5 &	 & 	 & 	16.2$^*$ & 	99.6 & 	7.8$^*$ & 	939.5 \\ 
case2736sp\_k\_api & 	 &	T.L. &	2.7 &	389.9 &	 & 	 & 	11.2$^*$ & 	118 & 	$\#$ & 	 \\ 
case2737sop\_k\_api & 	5.9 &	69.8 &	3.2 &	385.6 &	 & 	 & 	5.9 & 	94.1 & 	$\#$ & 	 \\ 
case2742\_goc\_api & 	18.6 &	1149.2 &	17.7 &	278.9 &	 & 	 & 	20.9$^*$ & 	700.8 & 	$\#$ & 	\\ 
case2853\_sdet\_api & 	\textbf{0.8} &	26.4 &	\textbf{0.4} &	76.6 &	 & 	 & 	MEM & 	 & 	 & 	 \\ 
case2869\_pegase\_api & 	\textbf{0.4} &	140.5 &	F &	 &	 & 	 & 	MEM & 	 & 	 & 	 \\ 
\hline
case3022\_goc\_api & 	 &	T.L. &	2.5 &	138 &	 & 	 & 		MEM & 	 & 	 & 	\\
case3120sp\_k\_api & 	 &	T.L. &	8.9 &	412.7 &	 & 	 & 		MEM & 	 & 	 & 	\\
case3375wp\_k\_api & 	3.8 &	316.1 &	 &	T.L. &	 & 	 & 	 	MEM & 	 & 	 & 	\\
case3970\_goc\_api & 	29.7 &	106.5 &	29.4 &	631.8 &	 & 	 & 		MEM & 	 & 	 & 	\\
case4020\_goc\_api & 	16.3 &	84.4 &	16.3 &	381.3 &	 & 	 & 	MEM	 & 	 & 	 & 	\\
case4601\_goc\_api & 	15.3 &	218.7 &	14.8 &	312.3 &	 & 	 & 	MEM & 	 & 	 & 	\\
case4619\_goc\_api & 	7 &	205.1 &	6.9 &	285.7 &	 & 	 & 	MEM	 & 	 & 	 & 	\\
case4661\_sdet\_api & 	1.6 &	162 &	1.4 &	455.7 &	 & 	 & 	MEM	 & 	 & 	 & 	\\
case4837\_goc\_api & 	6.6 &	185.7 &	6.5 &	984.9 &	 & 	 & 		MEM & 	 & 	 & 	\\
case4917\_goc\_api & 	5.2 &	1173.2 &	2.8 &	240.1 &	 & 	 & 		MEM & 	 & 	 & 	\\
case6470\_rte\_api & 	\textbf{0.3} &	201.5 &	 &	T.L. &	 & 	 & 		MEM & 	 & 	 & 	\\
case6495\_rte\_api & 	\textbf{0.8} &	2016 &	 & T.L.	 &	 & 	 & 	MEM	 & 	 & 	 & 	\\
case6515\_rte\_api & 	&	T.L.  &	 & T.L.	 &	 & 	 & 		MEM & 	 & 	 & 	\\
case9241\_pegase\_api & 	2.2 &	185.6 &	&	T.L.  &	 & 	 & 		MEM & 	 & 	 & 	\\
case9591\_goc\_api & 	13.7 &	87.6 &	13.7 &	1757.2 &	 & 	 &	MEM 	 & 	 & 	 & 	\\
case10000\_goc\_api & 	6.3 &	98.5 &	6.1 &	422.7 &	 & 	 & 	MEM	 & 	 & 	 & 	\\
case10480\_goc\_api & 	4.5 &	93 &	4.5 &	1255 &	 & 	 & 	MEM	 & 	 & 	 & 	\\
case13659\_pegase\_api & 	1.7 &	173.6 &	\textbf{0.9} &	1668.7 &	 & 	 & 	MEM	 & 	 & 	 & 	\\
case19402\_goc\_api & 	4.5 &	198.7 &	4.5 &	1768.1 &	 & 	 & 	MEM	 & 	 & 	 & 	\\
case24464\_goc\_api & 	3.8 &	183.6 &	3.8 &	2070.3 &	 & 	 & 	MEM	 & 	 & 	 & 	\\
case30000\_goc\_api & 	24.3 &	293 &	24.2 &	2797.1 &	 & 	 & 		MEM & 	 & 	 & 	\\
\hline
\hline
case3\_lmbd\_sad & 	1.9 &	0.0 &	0.1 &	0.1 &	  & 	  & 	3.6 & 	23.1 & 	0.0 & 	0.8 \\ 
case5\_pjm\_sad & 	0.0 &	0.0 &	0.0 &	0.1 &	  & 	  & 	0.0 & 	20.8 & 	 & 	 \\ 
case14\_ieee\_sad & 	0.3 &	0.0 &	0.3 &	0.1 &	  & 	  & 	0.1 & 	19.3 & 	 & 	 \\ 
case24\_ieee\_rts\_sad & 	4.4 &	0.0 &	4.4 &	0.1 &	 \textbf{0.2} & 	 1.9 & 	4.4 & 	18.8 & 	2.7 & 	2.1 \\ 
case30\_as\_sad & 	0.3 &	0.1 &	0.3 &	0.2 &	  & 	  & 	0.2 & 	20.4 & 	 & 	 \\ 
case30\_ieee\_sad & 	0.0 &	0.1 &	0.0 &	0.2 &	  & 	  & 	8.1 & 	28.7 & 	0.0 & 	1.7 \\ 
case73\_ieee\_rts\_sad & 	2.9 &	0.2 &	2.8 &	0.5 &	 \textbf{0.3} & 	16.8 & 	2.7 & 	19.8 & 	1.4$^*$ & 	7.4 \\ 
case118\_ieee\_sad & 	3.4 &	0.5 &	3.1 &	1.2 &	\textbf{0.0} & 	69.9 & 	3.1 & 	19.5 & 	1.9$^*$ & 	9.0 \\ 
case162\_ieee\_dtc\_sad & 	2.2 &	2.1 &	1.4 &	4.1 &	\textbf{0.0} & 	370.0  & 	5.4 & 	21.3 & 	1.2$^*$ & 	103.5 \\ 
case179\_goc\_sad & 	0.9 &	0.4 &	0.9 &	1.6 &	  & 	  & 	1.3 & 	19.5 & 	0.9$^*$ & 	11.9 \\ 
case240\_pserc\_sad & 	3.5 &	1 &	3.3 &	3.4 &	2.1 & 	T.L. & 	3.8$^*$ & 	20.7 & 	2.7$^*$ & 	32.9 \\ 
case300\_ieee\_sad & 	0.2 &	0.8 &	0.1 &	2.2 &	  & 	  & 	0.7$^*$ & 	22.4 & 	 & 	 \\ 
case500\_goc\_sad & 	5.6 &	4 &	5.6 &	4 &	 \textbf{0.4} & 	2398.7 & 	5.5$^*$ & 	29.2 & 	5.4$^*$ & 	87.6 \\ 
case588\_sdet\_sad & 	5.6 &	4.3 &	4.1 &	4.2 &	\textbf{ 0.4} & 	1859.7 & 	5.2$^*$ & 	24.7 & 	3.2$^*$ & 	63.5 \\ 
case793\_goc\_sad & 	4.9 &	4 &	4 &	9.5 &	\textbf{1.0} & 	7342.2 & 	5.1$^*$ & 	27.2 & 	3.0$^*$ & 	90.2 \\ 
case1354\_pegase\_sad & 	0.6 &	7.7 &	0.2 &	231.6 &	 & 	 & 	3.3 & 	37.6 & 	0.1 & 421.7	 \\ 
case1888\_rte\_sad & 	2.8 &	30.9 &	2.8 &	56.1 &	 & 	 & 	3.3$^*$ & 	50.9 & 	2.7$^*$ & 	547.5 \\ 
case2000\_goc\_sad & 	1.0 &	39.3 &	0.9 &	135.8 &	 & 	 & 	0.9$^*$ & 	73.8 & 	 & 	\\ 
case2312\_goc\_sad & 	 & T.L.	 &	1.6 &	109.1 &	 & 	 & 	3.3$^*$ & 	110.6 & 	1.2$^*$ & 	997.6 \\ 
case2383wp\_k\_sad & 	0.6 &	72.8 &	0.6 &	76.9 &	 & 	 & 	0.3$^*$ & 	91.3 & 	 & 	 \\ 
case2736sp\_k\_sad & 	0.4 &	1292.1 &	0.4 &	305.4 &	 & 	 & 	0.3$^*$ & 	100.2 & 	 & 	 \\ 
case2737sop\_k\_sad & 	 &	T.L. &	0.7 &	237.8 &	 & 	 & 	0.6$^*$ & 	91.0 & 	 & 	 \\ 
case2742\_goc\_sad & 	0.2 &	62.5 &	0.2 &	245.1 &	 & 	 & 	0.0$^*$ & 	580.3 & 	 & 	 \\ 
case2746wop\_k\_sad & 	 &	T.L. &	\textbf{0.7 }&	186.4 &	 & 	 & 	$\#$ & 	 & 	$\#$ & 	 \\ 
case2746wp\_k\_sad & 	 &	T.L. &	\textbf{0.5} &	139.3 &	 & 	 & 	$\#$ & 	 & $\#$	 & 	 \\ 
case2853\_sdet\_sad & 	1.1 &	29.9 &	\textbf{0.9} &	58.9 &	 & 	 & 	MEM & 	& 	 & 	 \\ 
case2869\_pegase\_sad & 	\textbf{0.1} &	39.3 &	F &	 &	 & 	 & 	MEM & 	& 	 & 	 \\ 
\hline
case3012wp\_k\_sad & 	\textbf{0.5} &	3571.9 &	\textbf{0.4} &	354.8 &	 & 	 & 		MEM & 	 & 	 & 	 \\ 
case3022\_goc\_sad & 	1.9 &	38.1 & \textbf{0.8} &	157.8 &	 & 	 & 		MEM & 	 & 	 & 	 \\ 
case3120sp\_k\_sad & 	 &	T.L. &	\textbf{0.6} &	159 &	 & 	 & 		MEM & 	 & 	 & 	 \\ 
case4020\_goc\_sad & 	8.2 &	263 &	8.2 &	231.7 &	 & 	 & 	MEM	 & 	 & 	 & 	 \\ 
case4601\_goc\_sad & 	6.1 &	436.8 &	6.1 &	221.3 &	 & 	 & 	MEM	 & 	 & 	 & 	 \\ 
case4619\_goc\_sad & 	1.7 &	77 &	1.7 &	935.4 &	 & 	 & 	MEM	 & 	 & 	 & 	 \\ 
case4661\_sdet\_sad & 	1.2 &	427.2 &	\textbf{0.9} &	456.3 &	 & 	 & 		MEM & 	 & 	 & 	 \\ 
case4917\_goc\_sad & 	1.6 &	107.3 &	\textbf{1.0} &	221.7 &	 & 	 & 	MEM	 & 	 & 	 & 	 \\ 
case6468\_rte\_sad & 	\textbf{0.6} &	191.5 &	F &	 &	 & 	 & 		MEM & 	 & 	 & 	 \\ 
case6470\_rte\_sad & 	\textbf{1.0} &	223 &	 &	T.L. &	 & 	 & 		MEM & 	 & 	 & 	 \\ 
case6495\_rte\_sad & 	13.3 &	500.3 &	F &	 &	 & 	 & 	MEM	 & 	 & 	 & 	 \\ 
case6515\_rte\_sad & 	7.1 &	318 &	F &	 &	 & 	 & 	 	MEM & 	 & 	 & 	 \\ 
case9241\_pegase\_sad & 	3.2 &	796.7 &	 &	T.L. &	 & 	 & 		MEM & 	 & 	 & 	 \\ 
case9591\_goc\_sad & 	9.5 &	179.6 &	9.5 &	388.5 &	 & 	 & 	MEM	 & 	 & 	 & 	 \\ 
case10000\_goc\_sad & 	6.3 &	582.2 &	6.3 &	282.3 &	 & 	 & 	MEM	 & 	 & 	 & 	 \\ 
case10480\_goc\_sad & 	1.1 & 3344.7	 &	1.1 &	1969.9 &	 & 	 & 	MEM	 & 	 & 	 & 	 \\ 
case13659\_pegase\_sad & 	2.1 &	168 &	 &	T.L. &	 & 	 & 	MEM	 & 	 & 	 & 	 \\ 
case19402\_goc\_sad & 	1.5 &	327.3 &	1.5 &	632 &	 & 	 & 		MEM & 	 & 	 & 	 \\ 
case24464\_goc\_sad & 	1.8 &	505.3 &	1.8 &	1213.3 &	 & 	 & 		MEM & 	 & 	 & 	 \\ 
case30000\_goc\_sad & 	10.9 &	379.4 &	10.9 &	755.4 &	 & 	 & 		MEM & 	 & 	 & 	 \\ 
\hline
\end{tabular}	
}
\vspace{0.05cm}
\caption{Optimality gaps from DSDP-RLT and Moment-Hierarchy}\label{Table:lb}
\end{table}	
}





\section{Conclusion}
In this paper, we present a comprehensive computational study on the ACOPF problem, benchmarking state-of-the-art open-source tools.
We also present methods and software solutions to compute lower bounds and optimality gaps on large-scale instances. The experiments show that using the general-purpose solver Ipopt along with the modeling language Gravity yields the best performance overall.

\section*{Acknowledgment}
We would like to thank Carleton Coffrin for providing us with scripts to run all the Julia experiments.
\bibliographystyle{IEEEtran}
\bibliography{ACOPF}
\end{document}